\documentclass[12pt]{article}
\usepackage{amssymb}
\usepackage{makeidx}
\usepackage{amsbsy,amsmath}
\usepackage{latexsym,amssymb,mathptm}
\usepackage{bm}

\newtheorem{Theorem}{Theorem}

\newtheorem{Proposition}[Theorem]{Proposition}

\title{{\bf On a Non-Classical Invariance Principle}}
\author{{Youri \sc Davydov${}^{1}$}\quad  and\quad { Vladimir \sc Rotar ${}^{2}$}}
\date
{\footnotesize  ${}^{1}$  Laboratoire Paul Painlev\'e - UMR 8524\\
 Universit\'e de Lille I - Bat. M2\\
 59655 Villeneuve d'Ascq, France \\
 Email: youri.davydov@univ-lille1.fr\\
\vspace{3pt}
${}^{2}$ \ Department of Mathematics and Statistics \\
 of the San Diego State University, USA and \\
 the Central Economics and Mathematics Institute \\
 of the Russian Academy of Sciences, RF  \\
 Email: vrotar@math.ucsd.edu }

\begin{document}

\maketitle

\begin{quote}
{\bf Abstract.}
We consider the invariance principle without the classical condition of
asymptotic negligibility of individual terms. More precisely, let r.v.'s $%
\{\xi _{nj}\}$ and $\{\eta _{nj}\}$ be such that
\[
E\{\xi _{nj}\}=E\{\eta _{nj}\}=0,\,\,E\{\xi _{nj}^{2}\}=E\{\eta
_{nj}^{2}\}=\sigma _{nj}^{2}, \,\,\sum_{j}\sigma _{nj}^{2}=1,
\]
and the r.v.'s $\{\eta _{nj}\}$ are normal. We set
\[
S_{kn} =\sum_{j=1}^{k}\xi _{nj},\,\,\,\,\,Y_{kn}=\sum_{j=1}^{k}\eta _{nj},\,
\,\,\, t_{kn} =\sum_{j=1}^{k}\sigma _{nj}^{2}.
\]
Let $X_{n}(t)$ and $Y_{n}(t)$ be continuous piecewise linear (or polygonal)
random functions with vertices at $(t_{kn},S_{kn})$ and $(t_{kn},Z_{kn})$,
respectively, and let $P_{n}$ and $Q_{n}$ be the respective distributions of
the processes $X_{n}(t)$ and $Y_{n}(t)$ in $\Bbb{C}[0,1]$.

The goal of the present paper is to establish necessary and sufficient
conditions for convergence of $P_{n} - Q_{n}$ to zero measure not involving
the condition of the asymptotic negligibility of the r.v.'s $\{\xi _{nj}\}$
and $\{\eta _{nj}\}$.


\end{quote}

\section{Introduction and results}

\renewcommand{\theequation}{\thesubsection.\arabic{equation}}

\subsection{Background and Motivation}

\setcounter{equation}{0}

The term ``non-classical'' concerns various limit theorems not
involving the condition of asymptotic negligibility of the
individual random variables (r.v.'s). To our knowledge, the
convergence of the distributions of sums of r.v.'s to the normal
distribution in the general situation, that is, without the
condition mentioned, was first considered by P. L\'{e}vy \cite{levy}
and M. Lo\'{e}ve \cite[Chapter VIII, Section 28]{loeve}. A developed
theory with necessary and sufficient conditions was built by V.M.
Zolotarev and his followers, V.M. Kruglov
and Yu.Yu. Machis; see, e.g., \cite{zolotarev2}, \cite{kruglov}, \cite{machis}%
, the monograph \cite{zolotarev}, the review part in \cite{r}, and
references therein. Note also that V.M. Kruglov considered the
Hilbert space case (see \cite{kruglov} and references in the papers
mentioned above.)

A somewhat different approach - see also comments below - that uses
different types of conditions, was suggested in \cite{r2} and
\cite{r}. In this paper, we proceed mainly from the framework of
\cite{r2} and \cite{r}.

In the case of normal convergence and finite variances, the simplest
result from \cite{r2} and \cite{r} may be stated as follows.

Let $\{\xi _{jn}\}$ be an array of independent r.v.'s such that
$E\{\xi
_{jn}\}=0,\,\,E\{\xi _{jn}^{2}\}=\sigma _{jn}^{2}<\infty $, and for each $n$%
,
\begin{equation}
\sum_{j}\sigma _{jn}^{2}=1.  \label{i1}
\end{equation}
Without loss of generality, we assume all $\sigma _{jn}\neq 0$.

Let $F_{jn}(x)$ be the distribution function (d.f.) of $\xi _{jn}$, and $%
\Phi _{jn}(x)$ be the normal d.f. with the same zero expectation and the
same variance; that is, $\Phi _{jn}(x)=\Phi (x/\sigma _{jn})$, where $\Phi
(x)$ is the standard normal d.f. \thinspace \thinspace Set $%
S_{n}=\sum_{j}\xi _{jn}$.

\begin{Proposition}
\label{pr1}(\cite{r2}) For
\begin{equation}
P(S_{n}\leq x)\rightarrow \Phi (x),\text{ for all }x\text{, as }n\rightarrow
\infty ,  \label{i2}
\end{equation}
it is necessary and sufficient that
\begin{equation}
\sum_{j}\int\limits_{||x|>\varepsilon }|x|\cdot |F_{jn}(x)-\Phi
_{jn}(x)|dx\rightarrow 0\text{, as }n\rightarrow \infty \text{, for any }%
\varepsilon >0.  \label{basic}
\end{equation}
\end{Proposition}

(This particular result is presented also in \cite{shiryaev} and \cite{r3}.)
It is easy to show (see, for example, \cite[p.310]{r3}) that in the
classical case where $\max_{j}\sigma _{jn}\rightarrow 0$, the Lindeberg
condition implies (\ref{basic}), so Lindeberg's theorem follows from
Proposition \ref{pr1}. On the other hand, condition (\ref{basic}) takes into
account possible proximity of the distributions of the r.v.'s to normal
ones. In particular, if $F_{jn}\equiv \Phi_{jn}$ and hence $P(S_n\leq
x)\equiv\Phi (x)$, then (\ref{basic}) becomes trivial.

It is worthwhile to note also that Proposition \ref{pr1} is
equivalent to Zolotarev's non-classical theorem from
\cite{zolotarev2} proved much earlier. In the framework of
\cite{zolotarev2}, the summands were directly divided into two
groups: those with ``small'' variances, and the rest. For the r.v.'s
from the former group, Lindeberg's condition was imposed, while the
summands from the latter group were required to be close to the
corresponding normal r.v.'s in L\'{e}vy's metric. Such a division
into two groups reflects the essence of the matter: ``small''
summands should be in the framework of the classical CLT, while
``large'' summands should be themselves close to normals. On the
other hand, condition (\ref{basic}) allows to treat the summands in
a unified way. Another difference between the theorem from
\cite{zolotarev2} and Proposition \ref{pr1} is that the latter uses
an integral metric.

In the sufficiency case, the result of Proposition \ref{pr1} was generalized
to the case of semi-martingales in Liptser and Shiryaev's paper \cite
{liptser}; see also Jacod and Shiryaev's book \cite[VII, 5b; VIII, 4c]{jacod}%
.

To generalize the result above to the case of convergence to distributions
different from normal, one may proceed as follows. Consider another array of
independent r.v.'s $\{\eta _{jn}\}$. We assume that for each $n$, the
numbers of terms for $\xi $'s and $\eta $'s in the arrays $\{\xi _{jn}\}$
and $\{\eta _{jn}\}$ are the same and, just for simplicity, are finite. Let $%
E\{\eta _{jn}\}=0,\,\,E\{\eta _{jn}^{2}\}=\sigma _{jn}^{2}$, and let $G_{jn}$
denote the distribution of $\eta _{jn}$. The problem is to establish
conditions under which
\begin{equation}
\prod_{j}F_{jn}-\prod_{j}G_{jn}\Rightarrow 0\text{ as }n\rightarrow \infty ,
\label{i3}
\end{equation}
where product of distributions is understood in the sense of convolution,
and convergence $\Rightarrow \,$itself is weak convergence (with respect of
all continuous bounded functions). At least formally, this is a more general
setup, since (\ref{i3}) does not presuppose the existence of limits for $%
\prod_{j}F_{jn}$ and $\prod_{j}G_{jn}$ separately. On the other hand, in the
particular case when $G_{jn}\equiv \Phi _{jn}$, (\ref{i3}) clearly coincides
with (\ref{i2}) in view of (\ref{i1}).

In the general situation (\ref{i3}), instead of (\ref{basic}), we consider
the condition
\begin{equation}
\sum_{j}\int\limits_{||x|>\varepsilon |}|x|\cdot
|F_{jn}(x)-G_{jn}(x)|dx\rightarrow 0,\text{ as }n\rightarrow \infty \text{,
for any }\varepsilon >0.  \label{basic1}
\end{equation}
In \cite{kir}, it was shown that when $G_{jn}$ are Poisson, (\ref{basic1})
remains to be a necessary and sufficient condition for the fulfillment of (%
\ref{i3}), however attempts to obtain a similar result in the general case
failed. The situation became clear when in \cite{rs1} and \cite{rs2} it was
proved that in general, relation (\ref{basic1}) is necessary for a more
stronger type of convergence. Namely, (\ref{basic1}) proves to be true if
and only if
\begin{equation}
\prod_{j\in B_{n}}F_{jn}-\prod_{j\in B_{n}}G_{jn}\Rightarrow 0\text{ as }%
n\rightarrow \infty ,  \label{con1}
\end{equation}
for any sequence $\{B_{n}\}$ of subsets of the indices $j$. See \cite{rs1}
and \cite{rs2} for detail; note also that in \cite{rs2} the case of infinite
variances is considered as well.

The fact that in the normal case, (\ref{i3}) and (\ref{con1}) occur to be
equivalent is connected with the fact that normal distributions are only
possible components of the decomposition of the normal law. The same
concerns the Poisson case, however in general, relations (\ref{i3}) and (\ref
{con1}) are certainly not equivalent.

Next, note that (\ref{con1}) deals with all possible partial sums,
so if we manage to establish the validity of this relation, it is
natural to continue and consider a more sophisticated problem,
namely, the asymptotic proximity of the distributions of the
partial-sum-processes based on the r.v.'s $\{\xi _{jn}\}$ and
$\{\eta _{jn}\}$.

The main goal of this note is to point out the fact that condition (\ref
{basic}) is \textit{necessary and sufficient } for the validity of
invariance principle in the case of Gaussian limiting processes in the
general, that is, non-classical setup. To our knowledge, this fact has not
been aired yet, though as we will see, in view of already known results, the
proof turns out to be not very difficult.

Note also that, as a matter of fact, we consider a slightly more general
problem of proximity of the distributions of the polygonal process generated
by the above r.v.'s $\xi _{jn}$ and the polygonal process generated by the
corresponding normal r.v.'s. In the classical case, when $\max_{j}\sigma
_{jn}\rightarrow 0$, such a result clearly corresponds to the classical
invariance principle of Donsker-Prokhorov (\cite{donsker}, \cite{prokhorov}%
), however without the condition mentioned we deal with a somewhat more
complicated situation.

We hope to consider a more general case of non-normal limiting distributions
in the next publication.

\subsection{Results}

\setcounter{equation}{0}

As was mentioned, we assume for simplicity that for each $n$, the
numbers of terms in each array, $\{\xi _{jn}\}$ or $\{\eta _{jn}\}$,
are finite. Suppose
all $\eta $'s are normal, so $G_{jn}(x)=\Phi _{jn}(x)=\Phi (x/\sigma _{jn})$%
. We again assume (\ref{i1}) to hold, and set

\begin{eqnarray}
S_{n} &=&\sum_{j}\xi _{jn},\,\,\,\,\,Y_{n}=\sum_{j}\eta _{jn},  \nonumber \\
S_{kn} &=&\sum_{j=1}^{k}\xi _{jn},\,\,\,\,\,Y_{kn}=\sum_{j=1}^{k}\eta _{jn},
\nonumber \\
t_{kn} &=&\sum_{j=1}^{k}\sigma _{jn}^{2}.  \label{n2}
\end{eqnarray}
Let $X_{n}(t)$ and $Y_{n}(t)$ be continuous piecewise linear (or polygonal)
random functions with vertices at $(t_{kn},S_{kn})$ and $(t_{kn},Y_{kn})$,
respectively. Let $\mathcal{P}_{n}$ and $\mathcal{Q}_{n}$ be the respective
distributions of the processes $X_{n}(t)$ and $Y_{n}(t)$ in $\Bbb{C=C}[0,1]$.

\begin{Theorem}
\label{th-main}Condition (\ref{basic}) is necessary and sufficient for
\begin{equation}
\mathcal{P}_{n}-\mathcal{Q}_{n}\Rightarrow 0  \label{cnv}
\end{equation}
(more precisely, to zero measure) weakly.
\end{Theorem}

Below, we show that the sequences $\{\mathcal{P}_{n}\}$ and $\{\mathcal{Q}%
_{n}\}$ are relatively compact, and hence in our case the above convergence
is equivalent to that in the L\'{e}vy-Prokhorov's metric $\pi $, that is, $%
\pi (\mathcal{P}_{n},\,\mathcal{Q}_{n})\rightarrow 0$. In general, when
compactness does not take place, and so to speak, ``parts of the
distributions move to infinity'', asymptotic proximity of distributions even
in the one-dimensional case may be defined in different ways, so the very
notion of proximity requires further analysis. We consider this question
separately in \cite{dr}.

We supplement Theorem \ref{th-main} by the following simple proposition. Let
for each $n$, the function $\sigma _{n}^{2}(t)=E\{X_{n}^{2}(t)\}$. Clearly, $%
\sigma _{n}(t)$ is continuous on $[0,1]$,

\[
\sigma _{n}^{2}(t_{kn})=\sum_{j=1}^{k}\sigma _{jn}^{2},
\]
and in each segment $[t_{(k-1)n},t_{kn}]$, the function $\sigma
_{n}^{2}(t)$ is a quadratic function.

\begin{Proposition}
\label{pr2}
 The process $Y_{n}(t)$
converges in distribution to a Gaussian process $Y(t)$ on $[0,1]$ such that $%
E\{Y(t)\}=0$ and $E\{Y^{2}(t)\}=\sigma ^{2}(t)$ if and only if
 for each $t\in [0,1]$,
\[
\sigma _{n}(t)\rightarrow \sigma (t).
\]

If $\max_{j}\sigma _{jn}\rightarrow 0$, then $\sigma ^{2}(t)=t$, and $Y(t)$
is the standard Wiener process. In general, the segment $[0,1]$ may be
divided into two sets, $A$ and $B$, with the following properties.

The set $A$ is a union of a finite or countable number of segments, and on
each such a segment the process $Y(t)$ is linear.

The set $B=[0,1]\diagdown A$, and if a segment $[a,b]\subset B$, then the
process $Y(a+s)-Y(a)$ is the standard Wiener process for $s\in [0,b-a]$.
\end{Proposition}

\section{Proofs}

The main issue is to prove the relative compactness of the measure
sequences $\{\mathcal{P}_{n}\}$ and $\{\mathcal{Q}_{n}\}$ (with
respect to weak convergence of distributions in $\Bbb{C}$). For
brevity, we omit sometimes the adjective ``relative''.

\subsection{Compactness in the normal case\label{cn}\label{cnc}}

\setcounter{equation}{0}

For the proof below, we need to consider a modification of the process $%
Y_{n}(t)$. For each $n=1,2,...$ , consider a partition of $[0,1)$ into
some intervals $[s_{(j-1)n},s_{jn})$ where $j=1,...,m_{n}\leq \infty $, and $%
0=s_{0n}<s_{1n}<...$. The number of intervals may be infinite,
points $s_{jn}$ may differ from the points $t_{jn}$ above.

Let $W_{n}(t)$ be a continuous piecewise linear process such that $%
W_{n}(0)=0 $, on each interval $[s_{(j-1)n},s_{jn})$ the trajectory of the
process is linear, and each increment $W_{n}(s_{jn})-W_{n}(s_{(j-1)n})$ is
either equal to zero, or to a normal r.v. $\zeta _{jn}$ with zero mean and a
variance of $s_{jn}-s_{(j-1)n}$. We prove the relative compactness of the
family of the distributions of \thinspace $W_{n}(t)$.

In accordance with a well known criterion (see, e.g., \cite{gs}), it
suffices to prove that

\begin{description}
\item  (A)\thinspace \thinspace \thinspace $\sup_{n}P\{|W_{n}(0)|>A\}%
\rightarrow 0$ as $A\rightarrow \infty $;

\item  (B)\thinspace \thinspace \thinspace there exist constants $a,b,c>0$
such that for any $n$ and $t,s\in [0,1]$,
\[
E\{|W_{n}(t)-W_{n}(s)|^{a}\}\leq c|t-s|^{1+b}.
\]
\end{description}

In our case, (A) is obvious. We verify (B) with $a=4,\,b=1$.

Set $v_{jn}^{2}=E\{(W_{n}(s_{jn})-W_{n}(s_{(j-1)n}))^{2}\}$. By the
definition of $W_{n}$, either $v_{jn}^{2}=0$, or $%
v_{jn}^{2}=s_{jn}-s_{(j-1)n}$.

If both points $t,s\in [s_{(k-1)n},s_{kn}]$ for some $k\geq 1$, and $%
v_{jn}^{2}\not{=}0$, then
\begin{equation}
E\{|W_{n}(t)-W_{n}(s)|^{4}\}=E\left\{ \left( \frac{|t-s|}{v_{kn}^{2}}\zeta
_{kn}\right) ^{4}\right\} \leq \frac{|t-s|^{4}}{v_{kn}^{8}}3v_{kn}^{4}\leq
3|t-s|^{2}  \label{n1}
\end{equation}
since in this case $|t-s|\leq v_{kn}^{2}$. On the other hand, if $v_{kn}=0$,
then $W_{n}(t)-W_{n}(s)=0$, and (\ref{n1}) is clearly true.

If $t=s_{kn}$ and $s=s_{mn}$ for some $k$ and $m>k$, then the r.v. $%
W_{n}(t)-W_{n}(s)$ is normal with a variance that does not exceed $%
s_{mn}-s_{kn}$. Then
\[
E\{|W_{n}(t)-W_{n}(s)|^{4}\}\leq 3\left( s_{mn}-s_{kn}\right)
^{2}=3|t-s|^{2}.
\]

In general, if $t\in [s_{(k-1)n},s_{kn}]$ and $s\in [s_{(m-1)n},s_{mn}]$ for
some $k$ and $m>k$, then in view of the above bounds,
\begin{eqnarray*}
E\{|W_{n}(t)-W_{n}(s)|^{4}\} \leq E\left\{ \left(
|W_{n}(t)-W_{n}(s_{kn})|+|W_{n}(s_{kn})-W_{n}(s_{(m-1)n})|\right. \right. \\
\left. \left. +|W_{n}(s_{(m-1)n})-W_{n}(s)|\right) ^{4}\right\} \\
\leq 27\left\{ E\left\{ |W_{n}(t)-W_{n}(s_{kn})|^{4}\right\} +E\left\{
|W_{n}(s_{kn})-W_{n}(s_{(m-1)n})|^{4}\right\} \right. \\
+\left. E\left\{ |W_{n}(st_{(m-1)n})-W_{n}(s)|^{4}\right\} \right\} \leq
243|t-s|^{2}.\,\,\blacksquare
\end{eqnarray*}

\subsection{Compactness of $\{\mathcal{P}_{n}\}$}

\setcounter{equation}{0}

First, note that in \cite[Lemma 2]{liptser}, relative compactness in the
non-classical situation was established in the general case of local \textit{%
martingales} with respect to weak convergence in $\Bbb{D}$. However, it is
not exactly what we need since we consider convergence in $\Bbb{C}$.

Certainly, once we consider continuous processes, and if limiting processes
are also continuous (which is true in our case), compactness in $\Bbb{D}$
implies convergence in $\Bbb{C}$. However, when considering piecewise linear
processes like $X_{n}(t)$ we loose the martingale property even when the
r.v.'s $\xi _{jn}$ are independent. On the other hand, if we switch to
piecewise constant processes, we have to consider convergence in $\Bbb{D}$,
which is not enough for us.

We believe that this is a technical obstacle and it may be somehow fixed,
but in any case, in our opinion, a self contained (and relatively short)
proof for the situation of independent summands would have an intrinsic
value. So, we provide this proof.

Thus, we establish relative compactness of $\{\mathcal{P}_{n}\}$ in $\Bbb{C}$
under condition (\ref{basic}).

Set $\mathbf{k}_{jn}=[t_{(j-1)n},\,t_{jn}]$, where the points $t_{jn}$ are
defined as in (\ref{n2}). For $\delta >0$, we define the process $%
X_{n}(t;\delta \,)$ as a result of replacement of the r.v.'s $\xi _{jn}$ by
the r.v.'s $\tilde{\xi}_{jn}=\xi _{jn}\mathbf{1}\{\sigma _{jn}^{2}>\delta \}$
in the definition of $X_{n}(t)$. (As usual, $\mathbf{1}\{A\}$ is the
indicator of a condition $A$.)

First, we show that for a fixed $\delta >0$, the family of the distributions
of $X_{n}(t;\delta \,)$ is compact. Indeed, denote by $\widetilde{\mathbf{k}}%
_{mn}^{\delta }=[r_{(m-1)n}^{\delta },r_{mn}^{\delta }]$ the segments $%
\mathbf{k}_{jn}$ where the process $X_{n}(t;\delta \,)$ is not constant. We
assume that $\widetilde{\mathbf{k}}_{mn}^{\delta }$ is on the left of $%
\widetilde{\mathbf{k}}_{(m+1)n}^{\delta }$. Since $\delta >0$, the number of
the segments $\widetilde{\mathbf{k}}_{mn}$ is finite. Denote this number by $%
q(n,\,\delta )$. Clearly, $q(n,\delta )\leq q=[1/\delta ]$ where $[a]$
stands for the integer part of $a$. It is convenient to think that always $%
m=1,...,q$, setting $\widetilde{\mathbf{k}}_{mn}^{\delta }=[1,1]$ for $%
m>q(n,\delta )$.

Clearly, there exists a subsequence $\widetilde{\mathbf{n}}=\{\widetilde{n}%
_{i}\}$ and segments $\widetilde{\mathbf{k}}_{m}^{\delta
}=[r_{(m-1)}^{\delta },r_{m}^{\delta }]$, $m=1,...,q$, such that
\[
\widetilde{\mathbf{k}}_{m\widetilde{n_{i}}}^{\delta }\rightarrow \,%
\widetilde{\mathbf{k}}_{m}^{\delta \,}\,\,\,\text{as\thinspace \thinspace
\thinspace \thinspace }i\rightarrow \infty ,
\]
(that is, the corresponding endpoints of the segments converge).

On the other hand, for each $\widetilde{\mathbf{k}}_{mn}^{\delta }$, the
distribution of the increment $X_{n}(r_{mn}^{\delta };\,\delta
)-X_{n}(r_{(m-1)n}^{\delta };\,\delta )$ is equal to a distribution $F_{jn}$
for some $j$. Then from the main condition (\ref{basic}) it follows that the
distribution of $X_{\widetilde{n}_{i}}(r_{m\widetilde{n}_{i}}^{\delta
};\,\delta )-X_{n}(r_{(m-1)\widetilde{n_{i}}}^{\delta };\,\delta )$ weakly
converges to the normal distribution with zero mean and the variance equal
to the length of $\widetilde{\mathbf{k}}_{m}^{\delta }$. (We skip a formal
proof of this fact. Because (\ref{basic}) is true for any $\varepsilon >0$,
we have convergence in the corresponding integral metric on any segments $%
[\varepsilon ,\infty )$ and $(-\infty ,-\varepsilon ]$. This implies weak
convergence. Since the limiting distribution is continuous, we have as a
matter of fact uniform convergence, but we do not need it.)

Since the distribution of the process $X_{n}^{\delta }(\cdot )$ is uniquely
specified by the finite dimensional distribution of the increments on the
segments $\widetilde{\mathbf{k}}_{m}^{\delta }$, we finally conclude that
the distribution of $X_{_{\widetilde{n}_{i}}}(\,\cdot \,;\,\delta )$ weakly
converges to the distribution of a continuous piecewise linear Gaussian
process $W(t;\,\delta )$ having points of growth only in the segments $%
\widetilde{\mathbf{k}}_{m}^{\delta }$ and such that the increments $%
W(r_{m}^{\delta };\,\delta )-W(r_{(m-1)}^{\delta };\,\delta )$ are normal
with zero mean and variance $r_{m}^{\delta }-r_{(m-1)}^{\delta }$.

Now, we proceed to a direct proof of compactness. Consider a sequence of
positive numbers $\delta _{n}\rightarrow 0$. As was shown, there exists a
subsequence $\mathbf{n}^{(1)}=\{n_{i}^{(1)}\}$ such that
\[
X_{n_{i}^{(1)}}(\,\cdot \,;\,\delta _{1})\stackrel{d}{\Rightarrow }W(\,\cdot
\,;\,\delta _{1})\,\,\text{as\thinspace \thinspace }i\rightarrow \infty ,
\]
where $\stackrel{d}{\Rightarrow }$ stands for weak convergence of the
corresponding distributions, and $W^{\delta _{1}}(\cdot )$ is a Gaussian
process of the type $W^{\delta }(\cdot )$ described above.

Similarly, we can choose a subsequence $\mathbf{n}^{(2)}$ of the sequence $%
\mathbf{n}^{(1)}$ such that
\[
X_{n_{i}^{(2)}}(\,\cdot \,;\,\delta _{2})\stackrel{d}{\Rightarrow }W(\,\cdot
\,;\,\delta _{2})\,\,\,\,\text{as\thinspace \thinspace }i\rightarrow \infty ,
\]
where $W^{\delta _{2}}(\cdot )$ is a Gaussian process with the same
properties as above. Continuing to reason in the same fashion, we
come to a nested sequence of subsequences $\mathbf{n}^{(1)}\supseteq
\mathbf{n}^{(2)}\supseteq ...\,\,\,$such that for all $k=1,2,...\,$,
\[
X_{n_{i}^{(k)}}(\,\cdot \,;\,\delta _{k})\stackrel{d}{\Rightarrow }W(\,\cdot
\,;\,\delta _{k})\,\,\,\,\text{as}\,\,\,\,i\rightarrow \infty .
\]

Next, consider the sequence of the Gaussian processes $\{W(\,\cdot
\,;\,\delta _{1}),\,\,W(\,\cdot \,;\,\delta _{2}),\,...\,\}$. By the result
of Section \ref{cn}, there exists a subsequence $m_{j}$ such that
\[
W(\,\cdot \,;\,\delta _{m_{j}})\stackrel{d}{\Rightarrow }W(\cdot )\,\text{,}
\]
where $W(\cdot )$ is a Gaussian process.

Now, we censor the sequence $\mathbf{n}^{(1)}\supseteq \mathbf{n}%
^{(2)}\supseteq ...\,\,\,$, choosing only $\mathbf{n}^{(m_{1})}\supseteq
\mathbf{n}^{(m_{2})}\supseteq ...\,\,\,$. By construction, we can choose a
sequence $n_{1},n_{2},...$ such that
\[
n_{1}\in \mathbf{n}^{(m_{1})},\,\,\,\,n_{2}\in \mathbf{n}^{(m_{2})},\,...%
\,,n_{i}\in \mathbf{n}^{(m_{i})},\,\,\,\,...\,\,
\]
and
\[
X_{n_{i}}(\,\cdot \,;\,\delta _{m_{i}})\stackrel{d}{\Rightarrow }W(\cdot
)\,\,\,\,\text{as}\,\,\,\,i\rightarrow \infty .
\]
$\,$ At the last step of the proof, we set $Z_{n}(t;\,\delta
)=X_{n}(t)-X_{n}(t;\,\delta )$, and consider the sequence of the processes $%
U_{i}(t)=Z_{n_{i}}(t;\,\delta _{m_{i}})$. Each process $U_{i}(t)$ \thinspace
\thinspace \thinspace \thinspace is a continuous process that is linear on
each segment $\mathbf{k}_{jn_{i}}$ and such that the variance of the
increment of the process on each $\mathbf{k}_{jn_{i}}$ does not exceed $%
\delta _{m_{i}}$. Since $\delta _{m_{i}}\rightarrow 0$ as $i\rightarrow
\infty $, all increments are asymptotically negligible. Formally, the
processes $\{U_{i}(t)\}$ are not exactly of the type appearing in the
classical invariance principle since for a finite number of segments $%
\mathbf{k}$ (with appropriate indices), the increments equals zero rather
than having a variance equal the length of $\mathbf{k}$. Nevertheless, the
proof of compactness may run exactly as, e.g., in the classical proof from
Prokhorov's paper \cite[Section 3.1]{prokhorov}.

Thus, the sequence of the distributions of $U_{i}(\,\cdot \,)$ is compact,
and so does the sequence of the distributions of $X_{n_{i}}(\,\cdot
\,;\,\delta _{m_{i}})$. It remains to observe that the processes $%
U_{i}(\,\cdot \,)$ and $X_{n_{i}}(\,\cdot \,;\,\delta _{m_{i}})$ are
independent.

\subsection{Proof of Theorem \ref{th-main}}

\setcounter{equation}{0}

\subsubsection{Necessity}

Let
\begin{equation}
\mathcal{P}_{n}-\mathcal{Q}_{n}\Rightarrow 0  \label{pro1}
\end{equation}
weakly in $\Bbb{C}[0,1]$. As was shown in Section \ref{cnc}, the sequence $\{%
\mathcal{Q}_{n}\}$ is compact. Then $\{\mathcal{P}_{n}\}$ is compact either.

Now, since $\sum_{j}\sigma _{jn}^{2}\equiv 1$, the marginal distribution
function for $Y_{n}$, i.e., $P(Y_{n}\leq x)\equiv \Phi (x)$. Hence, in view
of (\ref{pro1}),
\[
P(S_{n}\leq x)\rightarrow \Phi (x).
\]

By virtue of Proposition \ref{pr1}, this implies the validity of (\ref{basic}%
).

\subsubsection{Sufficiency}

Assume that condition (\ref{basic}) holds. Then, as was proved above, both
sequences, $\{\mathcal{P}_{n}\}$ and $\{\mathcal{Q}_{n}\}$, are compact.
Hence, it suffices to establish the convergence of the differences of all
finite-dimensional marginal distributions.

Let $t_{1}<t_{2}<...<t_{k}$ be points in $[0,1]$. Set $\mathbf{X}%
_{n}(t_{1},...,t_{k})=(X_{n}(t_{1}),...,X_{n}(t_{k}))$ and $\mathbf{Y}%
_{n}(t_{1},...,t_{k})=(Y_{n}(t_{1}),...,Y_{n}(t_{k}))$ and denote by $P_{n}%
\mathbf{(}t_{1},...,t_{k})$ and $Q_{n}\mathbf{(}t_{1},...,t_{k})$ the
distributions of the random vectors $\mathbf{X}_{n}(t_{1},...,t_{k})$ and $%
\mathbf{Y}_{n}(t_{1},...,t_{k})$, respectively. Both sequences, $%
\{P_{n}(t_{1},...,t_{k})$ and $Q_{n}(t_{1},...,t_{k})$, are compact.

We should prove that
\begin{equation}
P_{n}\mathbf{(}t_{1},...,t_{k})-Q_{n}\mathbf{(}t_{1},...,t_{k})\Rightarrow 0%
\text{.}  \label{pro2}
\end{equation}

Let the half interval $\mathbf{r}(j,n)=[t_{(j-1)n},t_{jn})$, and the
relations $t_{i}\in \mathbf{r}(m_{in},n)$, $i=1.,,,.k$, define the integers $%
m_{in}$. Then for $i=1,...,k$,
\begin{eqnarray}
X_{n}(t_{i}) &=&S_{(m_{in}-1)n}+\frac{t_{i}-t_{m_{i}n}}{\sigma _{m_{in}}^{2}}%
\xi _{m_{in}n},  \label{pro31} \\
Y_{n}(t_{i}) &=&Z_{(m_{in}-1)n}+\frac{t_{i}-t_{m_{i}n}}{\sigma _{m_{in}}^{2}}%
\eta _{m_{in}n}.  \label{pro32}
\end{eqnarray}
For each $n$, consider the random vectors
\begin{eqnarray}
&&\hspace{-45pt}\left( \sum_{j=1}^{m_{1n}-1}\xi _{jn},\,\,\,\xi
_{m_{1}n},\,\,\,\sum_{j=m_{1n}+1}^{m_{2n}-1}\xi _{jn},\,\,\,\xi _{m_{2}n}\,%
\mathbf{1}(m_{2}>m_{1}),\right.  \nonumber \\
&&\hspace{55pt}\left. ...,\sum_{j=m_{(k-1)n}+1}^{m_{kn}-1}\xi _{jn},\,\,\,\,\xi _{m_{k}n}%
\mathbf{1}(m_{k}>m_{k-1})\,\right) ,  \label{pro41}
\end{eqnarray}


and
\begin{eqnarray}
&&\hspace{-45pt}\left( \sum_{j=1}^{m_{1n}-1}\eta _{jn},\,\,\,\eta
_{m_{1}n},\,\,\,\sum_{j=m_{1n}+1}^{m_{2n}-1}\eta _{jn},\,\,\,\,\eta
_{m_{2}n}\,\mathbf{1}(m_{2}>m_{1}),\,...\,\right.  \nonumber \\
&&\hspace{55pt} \left. ...,\sum_{j=m_{(k-1)n}+1}^{m_{kn}-1}\eta _{jn},\,\,\,\,\eta
_{m_{k}n}\,\mathbf{1}(m_{k}>m_{k-1})\,\right)  \label{pro42}
\end{eqnarray}
where, by convention, $\sum_{a}^{b}=0$ for $a>b$.

Vectors (\ref{pro41}) and (\ref{pro42}) are those with independent
coordinates and are of the fixed dimension $2k$. Denote the $j$th
coordinates of these vectors by $\Psi _{jn}$, and $\Upsilon _{jn}$,
respectively, and set $\boldsymbol{\Psi }_{n}=(\Psi _{1n},...,\Psi _{2k,n})$%
, $\boldsymbol{\Upsilon }_{n}=(\Upsilon _{1n},...,\Upsilon _{2k,n})$. Let
the symbol $\mathcal{P}_{X}$ denote the distribution of a r.v. or a random
vector $X$.

First, note that the families of the distributions $\{\mathcal{P\,}_{%
\boldsymbol{\Psi }_{n}}\}$ and $\{\mathcal{P\,}_{\boldsymbol{\Upsilon }%
_{n}}\}$ are compact. Second, by results of \cite{rs1}-\cite{rs2} mentioned
in the Introduction, condition (\ref{basic}) implies that
\[
\prod_{j\in B_{n}}F_{jn}-\prod_{j\in B_{n}}\Phi _{jn}\Rightarrow 0
\]
weakly for any sequence $\{B_{n}\}$ of sets of indices. In particular, this
means that
\[
\mathcal{P}_{\Psi _{jn}}-\mathcal{P}_{\Upsilon _{jn}}\Rightarrow 0
\]
weakly for each $j$. Since the coordinates of the vectors $\boldsymbol{\Psi }%
_{n}$ and $\boldsymbol{\Upsilon }_{n}$ are independent, this implies that
\[
\mathcal{P\,}_{\boldsymbol{\Psi }_{n}}-\mathcal{P\,}_{\boldsymbol{\Upsilon }%
_{n}}\Rightarrow 0.
\]
On the other hand, in view of (\ref{pro31}) and (\ref{pro32}), each r.v.\ $%
X_{n}(t_{i})$ is a linear combination of the r.v.'s $\Psi _{jn}$, and each
r.v.\ $Y_{n}(t_{i})$ is the linear combination of the r.v.'s $\Upsilon _{jn}$
with \textit{the same coefficients }as for $X_{n}(t_{i})$. Together with the
compactness of $\mathcal{P\,}_{\boldsymbol{\Psi }_{n}}$ and $\mathcal{P\,}_{%
\boldsymbol{\Upsilon }_{n}}$, this leads to (\ref{pro2}). $\blacksquare
\vspace{0.1in}$

Since the sequence of the distributions $\{\mathcal{Q}_{n}\}$ is compact,
the proof of Proposition \ref{pr2} is straightforward, and we skip it.

\end{document}